\RequirePackage{fix-cm}
\documentclass[smallextended]{svjour3}       
\smartqed  
\usepackage{graphicx}
\usepackage{amsmath}
\usepackage{multirow}
\usepackage{mathtools}
\usepackage[english]{babel}
\usepackage{amsmath}
\usepackage{rotating}
\usepackage{array}
\usepackage{float}
\usepackage{subcaption}

\usepackage{caption}

\usepackage[numbers,compress]{natbib}

\begin{document}

\title{Rational Jacobi Collocation Method for the Solution of the Boundary Layer Flow of
an Eyring-Powell Non-Newtonian Fluid
over a Linear Stretching Sheet
}

\titlerunning{ Solving the Boundary Layer Flow of
an Eyring-Powell Non-Newtonian Fluid}        

\author{K. Parand*         \and
        S. Latifi\and
        M.M. Moayeri
}

\institute{K. Parand*(Corresponding Author) \at
              Department of Computer Sciences, Shahid Beheshti University, G.C. Tehran, Iran.\\
Department of Cognitive Modeling, Institute for Cognitive and Brain Sciences, Shahid Beheshti University, G.C. Tehran, Iran.\\
              \email{k\_parand@sbu.ac.ir}           
           \and
           S. Latifi \at
              Department of Computer Sciences, Shahid Beheshti University, G.C. Tehran, Iran.\email{s.latifi@mail.sbu.ac.ir} 
              \and
              M. M. Moayeri \at
               Department of Computer Sciences, Shahid Beheshti University, G.C. Tehran, Iran.               \email{mo.moayeri@mail.sbu.ac.ir} 
}
\date{Received: date / Accepted: date}
\maketitle

\label{intro}
\begin{abstract}

In this  paper, the  Rational Jacobi (RJ) collocation method is proposed to approximate the solution of the   boundary layer flow of an Eyring-Powell fluid over a stretching sheet. This equation is nonlinear  and by applying Quasilinearization method (QLM),  the equation is converted into a sequence of linear ordinary differential equations (ODE) converging to the solution of the nonlinear equation. Unlike other methods, instead of truncation in domain, the infinity condition is satisfied implicitly. As a result, using the proposed method,
the model is converted to a system of linear algebraic equations. The effect of different parameters on the velocity profile is also presented. 
\end{abstract}

\section{Introduction}

 In engineering processes, It happened to have stretching surface with a  velocity in an otherwise quiescent fluid. This stretching can cause a boundary layer viscous flow. This flow can be employed in metal extrusion, drawing of plastic films, polymer industry, glass blowing, paper production, wire drawing,  extrusion of polymer, artificial fiber, spinning of filaments, hot rolling, crystal growing and so forth.

In order to make different sheets, the percentage of melting is crucial and to have a special thickness and other desired formats, we need to consider stretching and cooling rate in this process.

In 1970, Crane \cite{crane1970flow} studied the flow of an incompressible viscous
fluid pasting a stretching sheet.

Mukhopadhyay \cite{mukhopadhyay2013slip} worked on the slip effects on Magnetohydrodynamic (MHD) boundary layer flow by an exponentially stretching sheet with suction/blowing and thermal radiation.  Bhattacharyya \cite{bhattacharyya2012effects} studied the effects of thermal radiation on micropolar fluid flow and heat transfer over a porous shrinking sheet.

Crane \cite{crane1970flow}, Gupta \cite{gupta1977heat}, Brady and Acrivos \cite{brady1981steady}, Wang and Usha \cite{wang1990liquid}, and Sridharan \cite{sridharan1995axisymmetric} all investigated the Flow of Newtonian fluids.   

Actually, Non-Newtonian fluids are  more important than Newtonian fluids and their mathematical models are more complex than Newtonian fluids.
The most frequent models that are used as  Non-Newtonian  models are power
law,  second grade, Maxwell and Oldroyd-B. 
Some uses of Non-Newtonian fluids are coal water, jellies, toothpaste,
ketchup, food products, inks, glues, soaps and they are also used in polymer solutions   \cite{hayat2013radiative}.
Unlike Newtonian fluids, mathematical systems for non-Newtonian fluids are of higher order and more difficult to deal with.

Nowadays, in addition to non-Newtonian,  Eyring-Powell model is another model of interest. The Eyring-Powell model has advantages over the other non-Newtonian fluid models because instead of the empirical relation, it deduced from Kinetic theory of liquids. Another advantage of it is that  it  reduces to Newtonian behavior for low and high shear rates \cite{rahimi2016solution}. In should be noted that in the calculation of fluid time scale at various polymer cases, Eyring-Powell  fluid models  are  more accurate and  consistent. Using Eyring-Powell model,  Patel and Timol \cite{patel2009numerical},  Hayat et al. \cite{hayat2012steady}, Rosca and Pop \cite{rocsca2014flow} and Jalil \cite{Jalil201373} investigated certain different facts and problems. 
Hayat et al \cite{hayaty} also observed the unsteady flow of Eyring-Powell fluid past an inclined stretching sheet. They write that unsteadiness in the flow is due to the time-dependence of the stretching velocity and wall temperature and in their study the corresponding  boundary layer equations are reduced into self-similar forms by suitable transformations then the analytic solutions are presented by homotopy analysis method (HAM). 
Furthermore,  Sugunamma et al. \cite{sugu} studied the effects of non uniform heat source/sink on unsteady flow of Powell-Eyring
fluid past an inclined stretching sheet with suction/injection effects ,and the corresponding equations are transformed into a system
of ordinary differential equations and solved it numerically using Runge- Kutta based shooting technique. \\

\subsection{Problems defined on Unbounded Domain}
 To solve the problems that are defined on unbounded intervals, numerical or semi-analytical methods can be applied. Among semi-analytical methods, Adomian Decomposition   \cite{tatari2007application}, Homotopy Perturbation \cite{shakeri2008solution}, Variational Iteration \cite{abdou2005variational} and Exp-function \cite{MaW2010X} methods can be pointed.  Finite Difference \cite{dehghan2008}, Finite Element  \cite{choi2016}, Meshfree \cite{Mirzaei2010} and  Spectral methods  \cite{parand2013,Rad2794} can be such categories of numerical methods.  
 
\subsection{Spectral Methods}
In recent years, the study of Spectral methods for solving ODEs and partial differential equations (PDEs) has attracted much attentions. The reason why these methods are used substantially is the accuracy that these methods provide in solving linear and nonlinear  phenomena like viscoelasticity, fluid mechanics, biology, physics and engineering. It is common to use orthogonal polynomials in the Spectral methods. There are so many authors being interested in utilizing these polynomials  \cite{
Bhrawyy2012,Doah2013,abd2013efficient,doha2014chebyshev,doha2015using}. In recent years, in addition to the Spectral methods and orthogonal  polynomials  and functions, other methods have also been applied \cite{parand2013solving,
parandhemami,ParandKHemamiM,
kazemRBF,parand2016operation,Rad2015363,Parand20144137}.\\
 Generally, in simple geometries for smooth problems, the Spectral methods offer exponential rates of convergence/Spectral accuracy. Nevertheless, Finite Element and Finite Difference methods just propose algebraic convergence rates. 
Galerkin, Collocation and Tau methods are the three most widely used approaches in  the Spectral methods. \\

To solve unbounded domain problems, several Spectral methods are used as different insights:\\

1) Direct approaches: In which orthogonal functions such as  Sinc, Hermite, Bessel and Laguerre are used \cite{shen2000,GuoBY2007,Guo2003,MaH2007}. The main feature of these functions is that they are defined over the unbounded intervals.\\

 2) Rational approximations: Where  some Spectral methods  on bounded domains are shifted into  unbounded domains \cite{BoydJP1987}. \\
 
 3) Mapping the  problem of interest: By which the problem defined on an unbounded domain is mapped into a bounded domain problem \cite{GuoBY2010}. \\
 
 4) Domain truncation: In this method, an enough large  number is assumed to be on behalf of infinity in the domain and all the conditions for infinity should be satisfied by this large number \cite{BoydJP1988}.\\
 
\subsection{Aim}
Due to the accuracy  needed to solve this problem defined on an unbounded domain,  the presented method is introduced. In this method the domain is not truncated and the infinity condition is satisfied implicitly. To increase the numerical solution accuracy, Spectral collocation method based on orthogonal RJ functions are used and RJ basis are also formed. These basis are made just by using an algebraic mapping to change the standard Jacobi polynomials defined on $[-1,1]$ into the RJ functions  on the domain $[0,\infty)$. 

The rest of this paper is structured as follows.
In the beginning,  certain basic facts about the definition and description  of the problem are reviewed, and the problem is formulated as a differential equation.  Then, 
RJ collocation method is described and the approximation is introduced as if it satisfies the infinity condition. For linearizing the nonlinear problem of the sort, we applied the QLM to ease the calculation and make the equation more straightforward. In the end, we solved the problem and reported the results and conclusion. 
\\

\section{Mathematical Formulation}
In this section, over a linear stretching sheet one  flow of a non-Newtonian Eyring-Powell fluid  is investigated. The linear velocity of stretched sheet 
is shown as $\tau X$, where $\tau$=b is the linear stretching velocity and $X$ is the distance from the slit \cite{rahimi2016solution}. In concise, we listed  parameters to use them in the next formulations.

\captionof{table}{Parameters and constants Description}
\begin{table}[H]\label{table:currentworktable}
\centering
 \begin{tabular}{|p{2cm}|p{4cm}|p{2cm}|p{3cm}|}
\hline
Parameter(s)   & Description & Parameter(s) & Description\\
\hline
$\xi_{ij}$ (kg/m $s^2$) &Shear-stress at the sheet&$\mu$ (kg/m s)& Dynamic viscosity of the fluid\\ 
\hline
$\gamma$&Material parameter&$H$&Material parameter\\
\hline
$u$, $v$ (m/s)&Non-dimensional velocity components along
$X$- and $Y$-axes&$X$,$Y$ (m)&Non-dimensional Cartesian coordinates\\
\hline
$\kappa$ (m2/s)& Kinematic viscosity of the fluid& $\rho$(kg/m3)& Fluid density\\
\hline
$\tau$ (m/s) &Stretching velocity&$\phi$&Stream function\\
\hline
$f(x)$&Dimensionless stream function&$x$&Similarity variable\\
\hline
$\delta$&Non-dimensional fluid parameter&$\varepsilon$&Non-dimensional fluid parameter\\
\hline

\end{tabular}
\end{table}

As in \cite{rahimi2016solution} mentioned, the extra stress tensor in an Eyring-Powell model is 
\begin{equation}\nonumber
\xi_{ij}=\mu\frac{\partial u_i}{\partial X_j}+\frac{1}{\gamma}sinh^{-1}\bigg(\frac{1}{H}\frac{\partial u_i}{\partial X_j}\bigg),
\end{equation}

and considering

\begin{equation}\nonumber
sinh^{-1}\bigg(\frac{1}{H}\frac{\partial u_i}{\partial X_j}\bigg)\cong \frac{1}{H}\frac{\partial u_i}{\partial X_j}-\frac{1}{6}\bigg(\frac{1}{H}\frac{\partial u_i}{\partial X_j}\bigg)^3 ,\bigg|\frac{1}{H}\frac{\partial u_i}{\partial X_j}\bigg| <<1.
\end{equation}

For the fluid based on Eyring-Powell model, the
equation of continuity and the $X$-momentum equation is 
simplified as
\begin{equation}\nonumber
\frac{\partial u}{\partial X}+\frac{\partial v}{\partial Y}=0,
\end{equation}

\begin{equation}\label{equation:e1}
u\frac{\partial u}{\partial X}+v\frac{\partial v}{\partial Y}=\frac{\partial^2 u}{\partial Y^2}\bigg(\kappa+\frac{1}{\rho\gamma H}\bigg)-\frac{1}{2\rho\gamma H^3}\bigg(\frac{\partial u}{\partial Y}\bigg)^2\frac{\partial^2 u}{\partial Y^2},
\end{equation}

and the boundary conditions are
\[   
u = 
     \begin{cases}
       bX,&Y=0, \\
      0, & Y \rightarrow \infty,\\
     \end{cases}
\]
additionally, if we set
$$u=\frac{\partial \phi}{\partial Y},v=-\frac{\partial \phi}{\partial X},$$

\begin{equation}\label{equation:e2}
\phi=Xf(x)(b\kappa)^{0.5},x=Y(\frac{b}{\kappa})^{0.5}.
\end{equation}
where $x$ is the similarity variable. Consequently, Eq.(\ref{equation:e1}) and Eq.(\ref{equation:e2}) lead to
\begin{equation}\label{equation:eq3}
f(x)f''(x)+(1+\varepsilon)f'''(x)-\varepsilon\delta f''(x)^2f'''(x)-f'(x)^2=0,
\end{equation}
and boundary conditions will be enhanced as $f(0)=0$ and 
\[
    f'(x)=\begin{cases}
       1,& x \rightarrow 0, \\
      0, & x \rightarrow \infty.\\
     \end{cases}
\]
The constant parameters $\delta$, $\varepsilon$ are  the Non-dimensional fluid parameter which are considered constant in numerical examples. In fact, these quantities have the following definitions \cite{rahimi2016solution}

\begin{equation}
\varepsilon=\frac{1}{\mu\gamma H},~~\delta=\frac{b^3X^2}{2H^2\kappa}.
\end{equation}
 
\section{Preliminaries}
\subsection{Jacobi Polynomials}
 Because of the features that Jacobi polynomials provide, these polynomials are considered to solve different types of problems \cite{ParandLatifi,bhrawyzaky2015,bhrawyAlzaidy,BhrawyDoha2015,BhrawyAHDohaEH2015,Doha2015pseudo}.
The standard Jacobi polynomial of degree $n$ with $\alpha,\beta>-1$ is as follows \cite{shen2011spectral}
$$J_{n}^{\alpha,\beta}(x)=\bigg(\frac{\Gamma(n+\alpha+1)}{n!\Gamma(n+\alpha+\beta+1)}\bigg)\sum_{i=0}^{n}\binom{n}{i}\frac{\Gamma(n+i+\alpha+\beta+1)}{\Gamma(i+\alpha+1)}{\frac{(x-1)}{2}}^i,$$\\
with the properties as:
\begin{equation}\nonumber
J_{n}^{\alpha,\beta}(-x)=(-1)^nJ_{n}^{\beta,\alpha}(x), J_{n}^{\alpha,\beta}(-1)=\frac{(-1)^n\Gamma(n+\beta+1)}{n!\Gamma(\beta+1)}, J_{n}^{\alpha,\beta}(1)=\frac{\Gamma(n+\alpha+1)}{n!\Gamma(\alpha+1)}.
\end{equation}

and its weight function is  $w^{\alpha,\beta}(x)=(1-x)^{\alpha}(1+x)^{\beta}$.
\\These polynomials are orthogonal on $[-1,1]$
\begin{equation}\nonumber
\int_{-1}^{1}{J_{n}^{\alpha,\beta}(x)J_{m}^{\alpha,\beta}(x)}w^{\alpha,\beta}(x)=\delta_{m,n}\frac{2^{\alpha+\beta+1}\Gamma(n+\alpha+1)\Gamma(n+\beta+1)}{(2n+\alpha+\beta+1)\Gamma(n+1)\Gamma(n+\alpha+\beta+1)},
\end{equation}
where $\delta_{m,n}$ is 
\[
    \delta_{m,n}=\begin{cases}
       1,& m=n, \\
      0, & m \neq n.\\
     \end{cases}
\]
The set of Jacobi polynomials makes a complete $L_{w^{\alpha,\beta}}^2[-1,1]$ orthogonal system.
\subsection{ Rational  Jacobi (RJ) Functions} 
 In spite of $[-1,1]$, the equation of interest is defined over an unbounded domain; To have all the properties of Jacobi polynomials on unbounded domain, we can use the mapping $\frac{x-L}{x+L}$ \cite{ParandKMazaherithomas}. Using this mapping will lead to creation of RJ  Basis defined as 
 
 $$RJ_{n,L}^{\alpha,\beta}(x)=J_{n}^{\alpha,\beta}(\frac{x-L}{x+L}).$$
These new RJ functions, inherit most of the useful characteristics of the standard Jacobi polynomials. Of these characteristics can mention: simplicity, exponential or infinite order convergence, completeness and having  real roots that are spread on $[0,\infty)$.
 By this, the domain will be $[0,\infty)$. Putting this mapping relation in standard Jacobi polynomials formula, we have

$$RJ_{n,L}^{\alpha,\beta}(x)=\bigg(\frac{\Gamma(n+\alpha+1)}{n!\Gamma(n+\alpha+\beta+1)}\bigg)\sum_{i=0}^{n}\binom{n}{i}\frac{\Gamma(n+i+\alpha+\beta+1)}{2\Gamma(i+\alpha+1)}{{\bigg((\frac{x-L}{x+L})-1\bigg)}}^i,$$\\
     and
$$ RJ_{n,L}^{\alpha,\beta}(-x)=(-1)^nRJ_{n,L}^{\beta,\alpha}(x), RJ_{n,L}^{\alpha,\beta}(0)=\frac{(-1)^n\Gamma(n+\beta+1)}{n!\Gamma(\beta+1)}, RJ_{n,L}^{\alpha,\beta}(\infty)=\frac{\Gamma(n+\alpha+1)}{n!\Gamma(\alpha+1)}.$$
The above mentioned weight function  for these new rational functions is  $w_{L}^{\alpha,\beta}(x)=(1-(\frac{x-L}{x+L}))^{\alpha}(1+(\frac{x-L}{x+L}))^{\beta}\frac{2L}{(x+L)^2}$. This weight function can be simplified to the form \\
$$w_{L}^{\alpha,\beta}(x)=\frac{2^{\alpha+\beta+1}x^{\beta}L^{\alpha+1}}{(x+L)^{\alpha+\beta+2}}.$$
Obviously, this weight function is  non-negative, integrable and in fact, that is a real-valued weight function over $[0,\infty)$.\\Moreover, RJ functions are orthogonal on $[0,\infty)$
\begin{equation}\nonumber
\int_{0}^{\infty}{RJ_{n,L}^{\alpha,\beta}(x)RJ_{m,L}^{\alpha,\beta}(x)}w_{L}^{\alpha,\beta}(x)=\delta_{m,n}\frac{2^{\alpha+\beta+1}\Gamma(n+\alpha+1)\Gamma(n+\beta+1)}{(2n+\alpha+\beta+1)\Gamma(n+1)\Gamma(n+\alpha+\beta+1)}.
\end{equation}
The set of RJ functions makes a complete $L_{w_L^{\alpha,\beta}}^2[0,\infty)$ orthogonal system and for any $y(x)\in L_{w_L^{\alpha,\beta}}^2[0,\infty),$ there is an expansion as follows.\\
$$y(x)=\sum_{j=0}^{\infty} b_jRJ_{j,L}^{\alpha,\beta}(x).$$

\section{Quasilinearization Method (QLM)}\label{section:QLM}
In reality, many of equations are nonlinear and difficult to deal with. One way to make these equations linear and utilize their linearity property is to use QLM.
QLM is based on Newton-Raphson method  \cite{khan2006generalized,el2006generalized,vatsala1999generalized,el2007quasilinearization}
 by which in an iterative view, we suppose that the solution of $i^{th}$ step is known and we attempt to find the solution in $(i+1)^{th}$ step.  Although it is an iterative method, the QLM is not perturbative and  the fast quadratic convergence is the most dominant quality of it  \cite{MandelzweigVB}.

This method was developed  by Bellman and Kalaba \cite{KalabaR,BellmanRE} and can be applied to solve  miscellaneous nonlinear ODEs or partial differential equations (PDEs)   in such different fields of sciences  like physics,  engineering and  biology.  As a matter of fact, by using QLM, a nonlinear equation can be converted into a linear equation in which the solution  converges to the solution of the nonlinear equation. It is worth mentioning that in  \cite{delkhoshtomas}, the convergence of QLM has been discussed and proved.
The  QLM  iteration  requires  an  initialization  or  "initial  guess" that  is  chosen  from  physical  and  mathematical  considerations  or  the  boundary
conditions in the equation.

To know how it works, imagine that we have one particular differential equation as:
 \begin{equation}\nonumber
 L^{(m)}(y(x))=h(x,y(x),y^{(1)}(x),y^{(2)}(x),...,y^{(m-1)}(x)).
 \end{equation}
where $L^{(m)}$ is the differential operator of order $m$, and $h(x,y(x),y^{(1)}(x),y^{(2)}(x),...,y^{(m-1)}(x))$  is a nonlinear function containing $x,y(x),y'(x),...,y^{(m-1)}(x)$ .\\  
 The linear equation after applying QLM will be

\begin{eqnarray}
&L^{(m)}(y_{i+1}(x))=h\big(x,y_{i}(x),y_{i}^{(1)}(x),y_{i}^{(2)}(x),...,y_{i}^{(m-1)}(x)\big)+\nonumber\\
&\sum_{j=0}^{m-1}{\big(y_{i+1}^{(j)}(x)-y_{i}^{(j)}(x)\big)h_{y^{(j)}}\big(x,y_{i}(x),y_{i}^{(1)}
(x),y_{i}^{(2)}(x),...,y_{i}^{(m-1)}(x)\big).}\nonumber
 \end{eqnarray}
  Convergence of the QLM proved by Bellman \& Kalaba \cite{BellmanRE}  and Mandelzweig  \& Tabakin \cite{mandelzweig2001quasilinearization}.
\\If  $\delta y_{i+1}(x)\equiv y_{i+1}(x)- y_i(x)$
be the difference between two subsequent iterations, then, 
$\| \delta y{i+1}\| \leq  k\| \delta y_i \|^2$, where k is a real constant. Therefore the convergence rate is of order 2 i.e. $O(h^2)$. It worth mentioning that the order of convergence in numerical methods can
be affected and reduced by the errors in computational methods, machine errors, roundoff errors,
and so forth.\\
Additionally we have the relationship\\
$$\| \delta y_{i+1} \| \leq ( \| \delta y_1 \| )^{(2^n)} /k $$
and we can say that the convergence depends on $q_1=k\|\delta y_1 \|= k \| \delta y_1-y_0 \|$.  For convergence it is sufficient that just one of $q_m=k\|\delta y_m\|$ be small enough. Even if the first convergent coefficient $q_1$ is large, a well chosen value for $y_0$ results in a small value for at least one of $q_m,m\geq 2$ and this also reduces the order of convergence. \cite{BellmanRE,mandelzweig2001quasilinearization}.
\\

  Considering Eq.(\ref{equation:eq3}),  the Eyring-Powell equation  of interest is  linearized and Eq.(\ref{equation:eq4}) is achieved. For the purpose of simplicity  [(.) $\equiv$ $\big(x,f_i(x),f'_i(x),f''_i(x)\big)$] is assumed (See Appendix).\\  
  \begin{equation}\label{equation:eq4}
  g_1(.)+g_2(.)f_{i+1}(x)+g_3(.)f'_{i+1}(x)+g_4(.)f''_{i+1}(x)-f'''_{i+1}(x)=0,
  \end{equation}
  where $g_j(.),j=1,2,3,4$ are all known and $f_{i+1}(x)$ is unknown and at the initial step $f_0:=1$. Clearly, the Eq.(\ref{equation:eq4}) is linear. From now on, instead of the former nonlinear equation, this linear equation is under consideration and in the next sections we attempt to solve it numerically.

 \section{Rational  Jacobi (RJ) Collocation Method}
In collocation method, It matters that how collocation points are chosen. Choosing these points can even affect convergence and efficiency.
Here in this paper, we approximate the solution  in an expansion like Eq.(\ref{equation:eq5}).
If $y(x)$ be the solution of our ODE, by continuing in $N+1$ points, its approximated expansion in iteration $(i+1)^{th}$ i.e. $f_{i+1}(x)$ will be 
  \begin{equation}\label{equation:eq5}
y(x)\cong f_{i+1}(x)=  \frac{x}{x^2+1}+\frac{x^2}{x^2+1}\sum_{j=0}^{N}{a_j^{i+1}RJ_{j,L}^{\alpha,\beta}(x)},
  \end{equation}
  in which all the following conditions are satisfied. 
  $$\lim_{x\to 0}f_{i+1}(x)=0
,\lim_{x\to0}f_{i+1}'(x)=1
    ,\lim_{x\to\infty}f_{i+1}'(x)=0.$$
Considering Eq.(\ref{equation:eq4}) we can show the residual function as (in order to read more  on this residual function see Appendix)
\begin{equation*}
Res(x;a_0^{i+1},a_1^{i+1},...,a_{N}^{i+1})=-f'''_{i+1}(x)+g_1(.)+
g_2(.)f_{i+1}(x)+g_3(.)f'_{i+1}(x)
+g_4(.)f''_{i+1}(x),
\end{equation*}
In which, $g_k(.), k=1,2,3,4$ are known (defined in Appendix) and the $f_{i+1}(x)$ is unknown and we are to approximate it. To find the unknown $a_j^{i+1},j=0,..,N$ we appoint the RJ collocation points($x_k$, $k=0,...,N$) and solve the system in Eq.(\ref{equation:eq6}). These collocation points are the roots of $RJ_{N+1,L}^{\alpha,\beta}(x)$.  
\begin{equation}\label{equation:eq6}
Res(x_k;a_0^{i+1},a_1^{i+1},...,a_{N}^{i+1})=0,k=0,...,N.
\end{equation}

\section{Numerical Results}
\begin{figure}[H]\label{fig:diff N}
    \centering
    \begin{subfigure}[b]{0.3\textwidth}
        \centering
        \includegraphics[width=\textwidth]{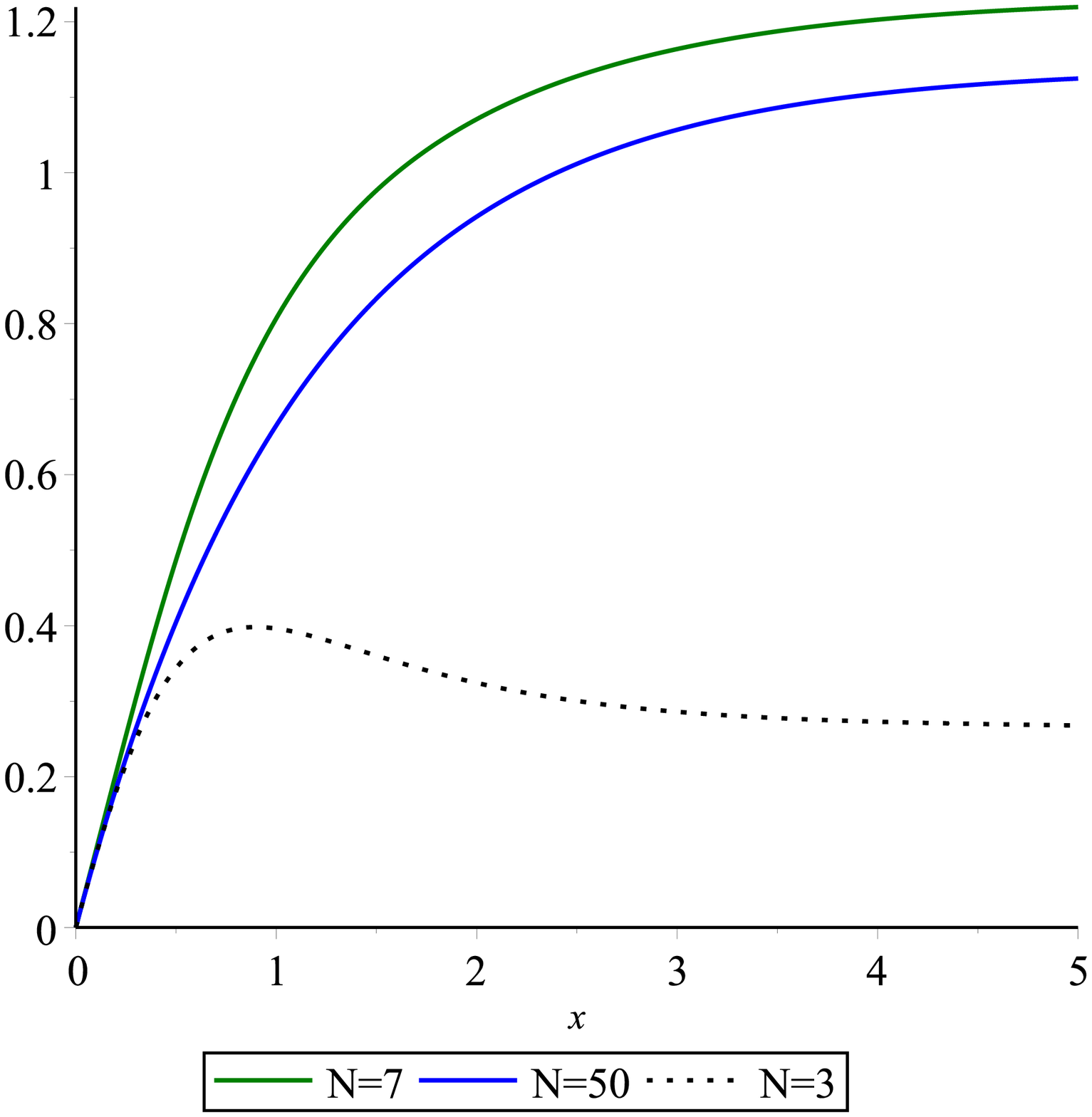}
        \caption{$f(x)$}
        \label{fig:a}
    \end{subfigure}
    \hfill
    \begin{subfigure}[b]{0.3\textwidth}
        \centering
        \includegraphics[width=\textwidth]{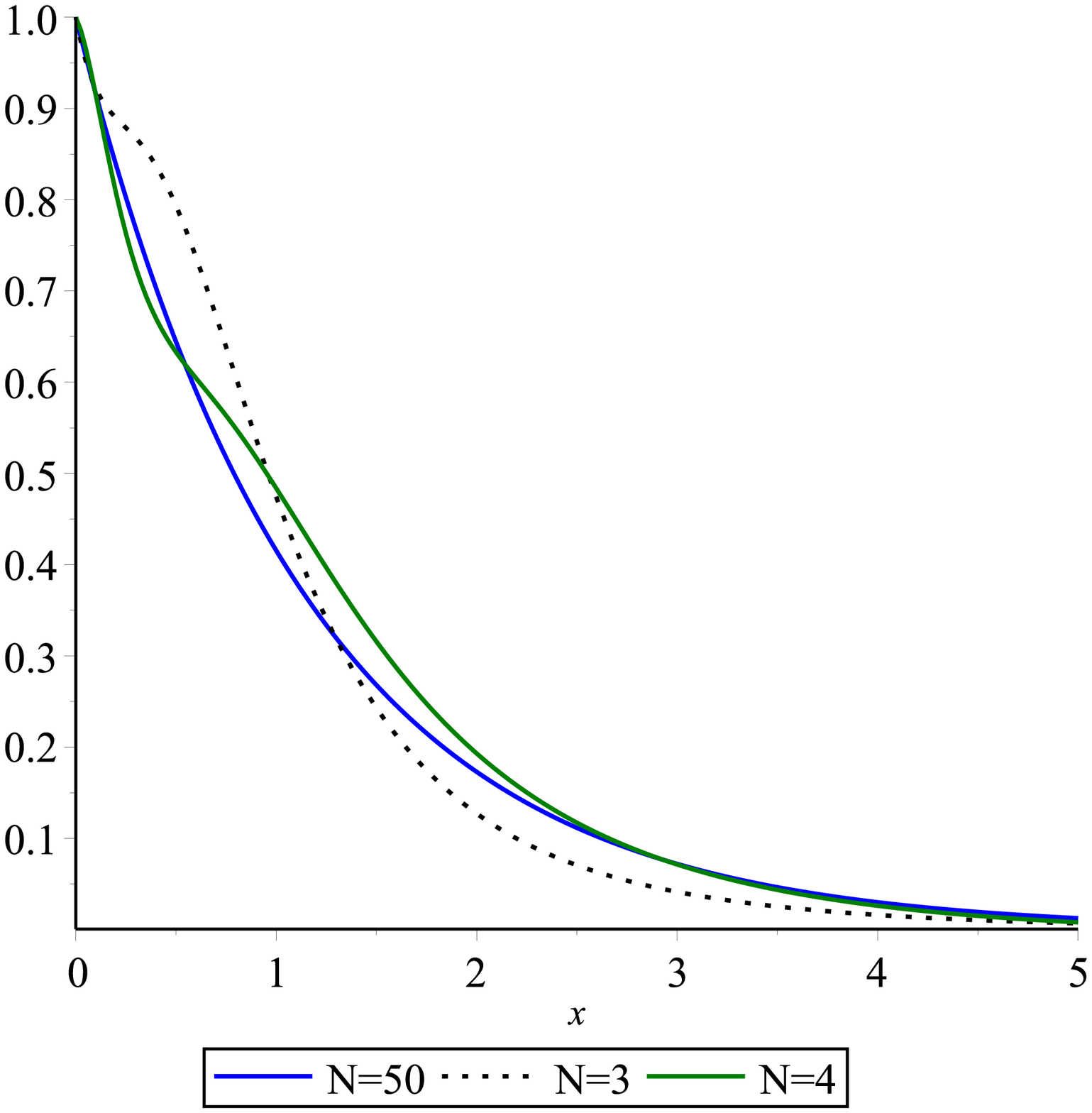}
        \caption{$f'(x)$}
        \label{fig:b}
    \end{subfigure}
    \hfill
    \begin{subfigure}[b]{0.3\textwidth}
        \centering
        \includegraphics[width=\textwidth]{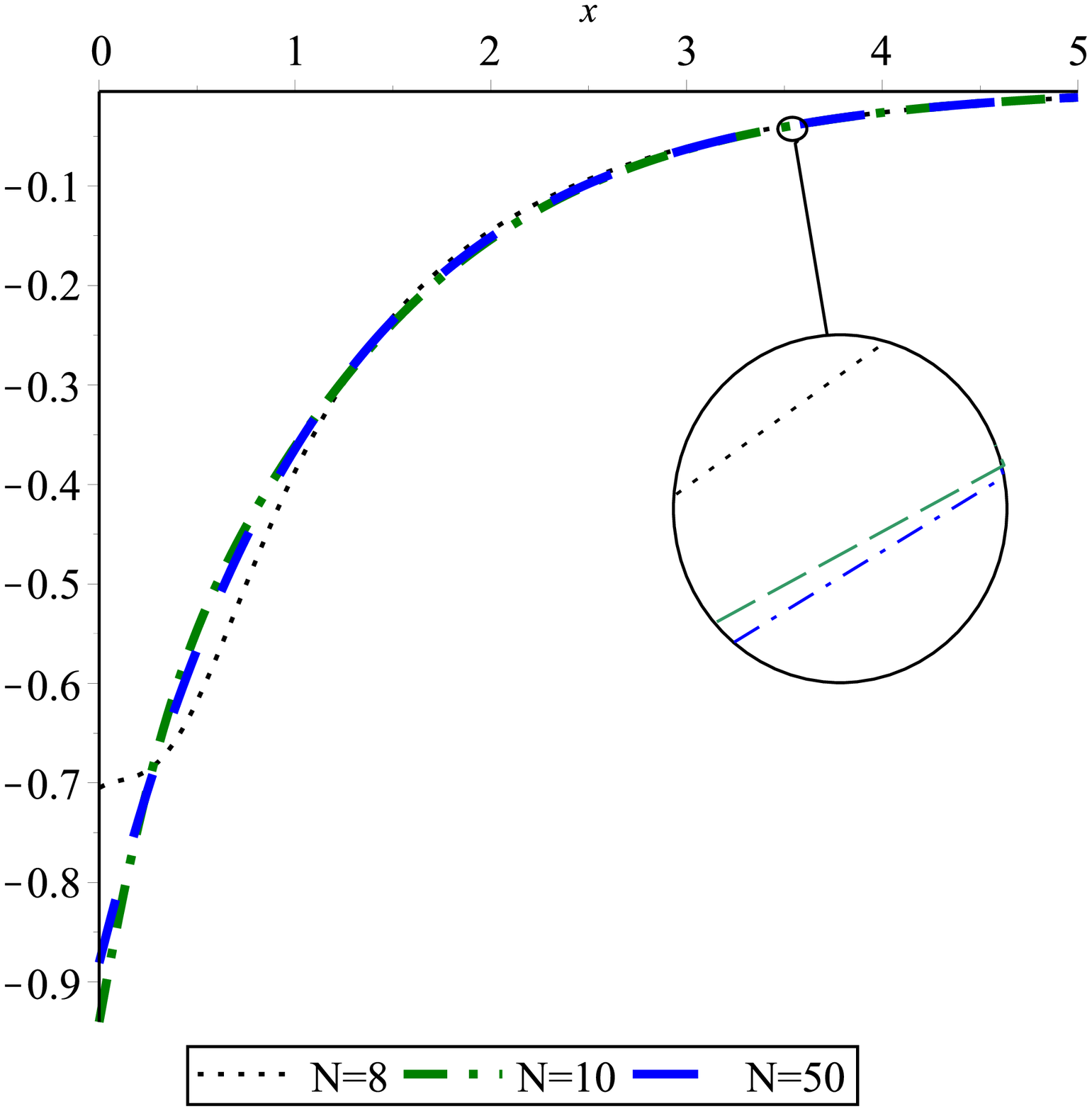}
        \caption{$f''(x)$}
        \label{fig:c}
    \end{subfigure}
     \hfill
    \begin{subfigure}[b]{0.4\textwidth}
        \centering
        \includegraphics[width=\textwidth]{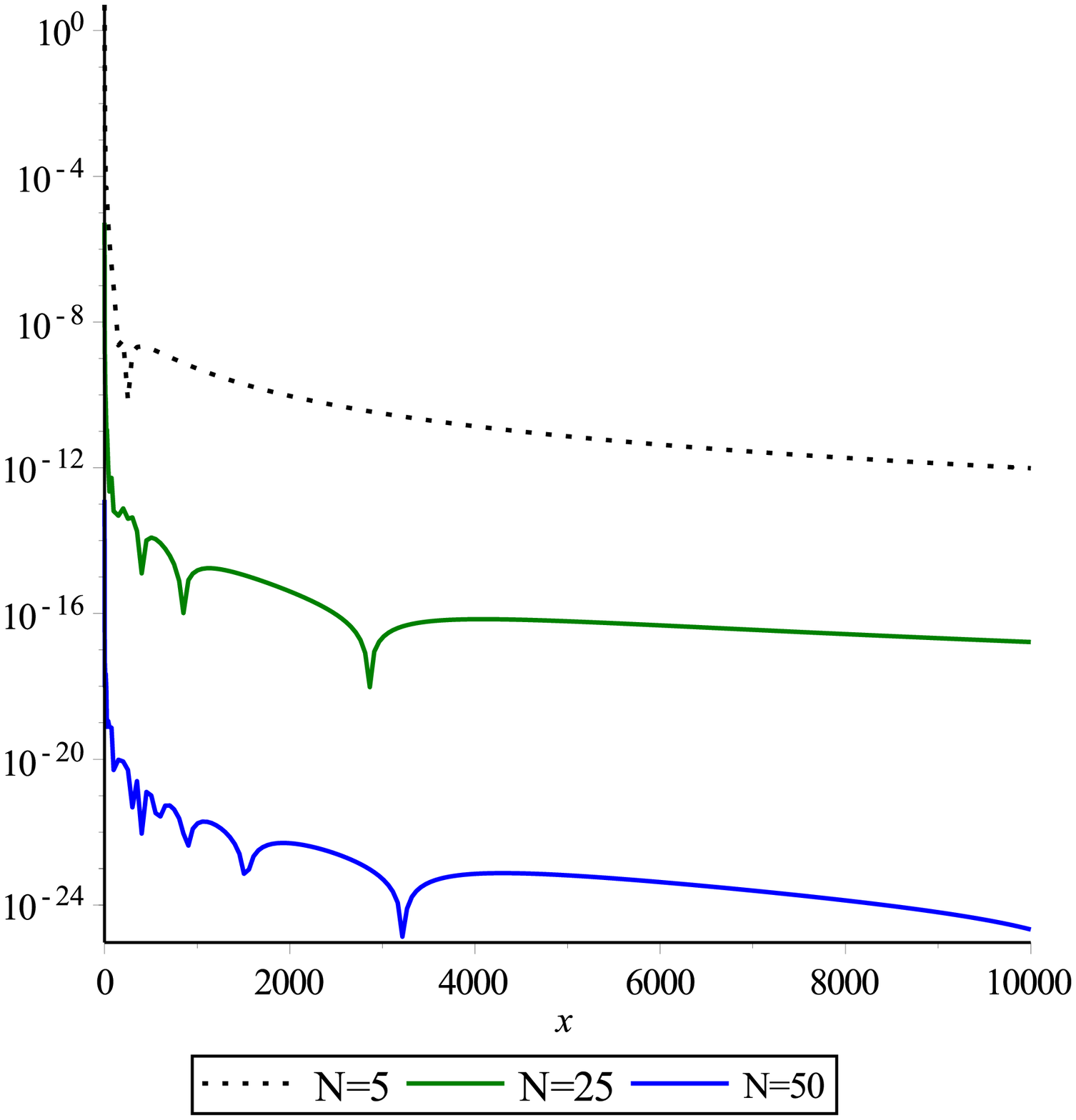}
        \caption{ Logarithm of the absolute value of residual.}
        \label{fig:d}
    \end{subfigure}
    \caption{ Diagram of $f(x)$, $f'(x)$, $f''(x)$ and $Log(|res(x)|)$ in the $25$th iteration for different $N$. $\delta=0.1$, $\varepsilon=0.3$, $\alpha=1,\beta=1$, $L=15$. 
}

\end{figure}

\begin{figure}[H]\label{fig: diff iter}
    \centering
    \begin{subfigure}[b]{0.3\textwidth}
        \centering
        \includegraphics[width=\textwidth]{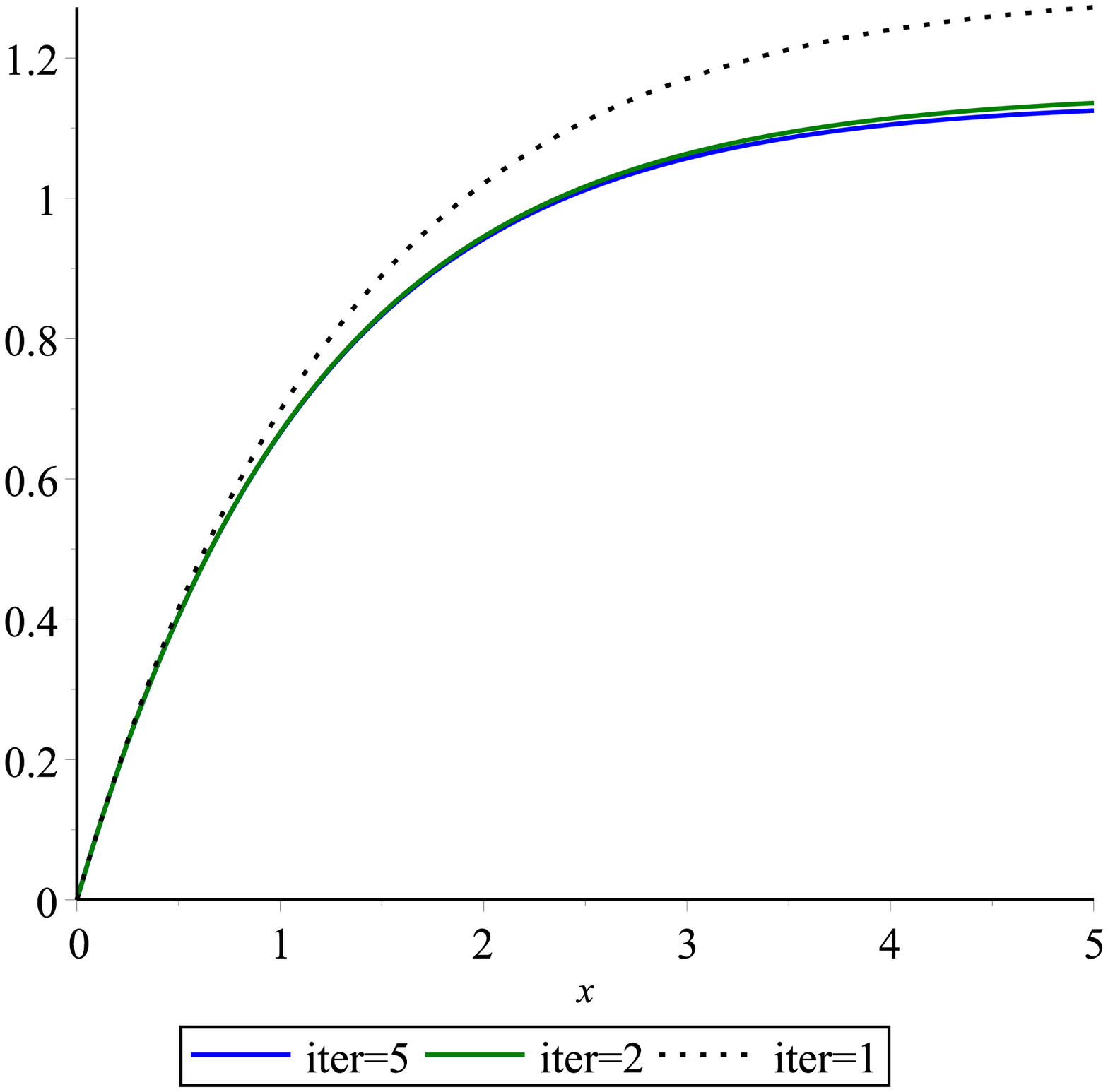}
        \caption{$f(x)$}
        \label{fig:e}
    \end{subfigure}
    \hfill
    \begin{subfigure}[b]{0.3\textwidth}
        \centering
        \includegraphics[width=\textwidth]{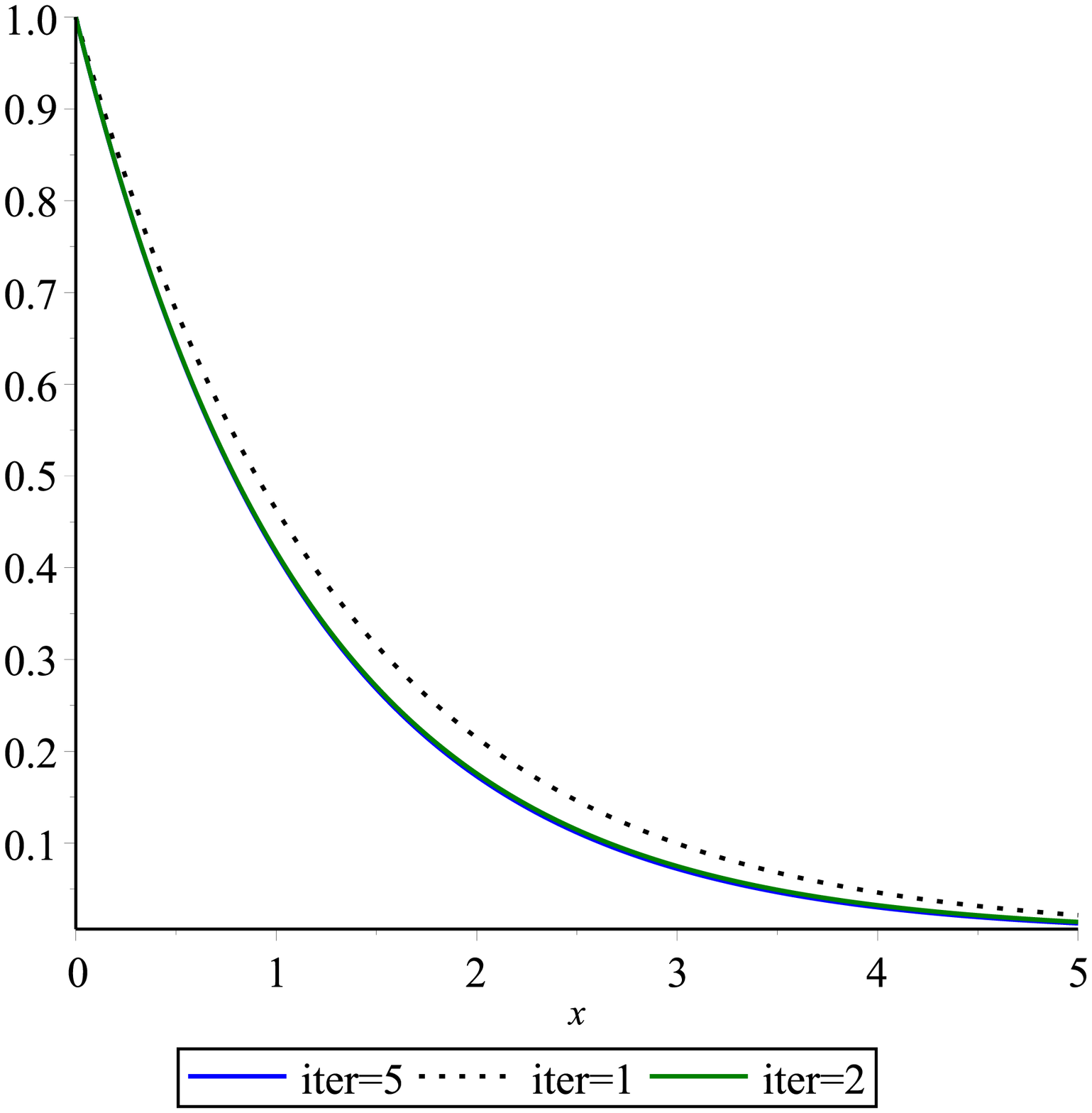}
        \caption{$f'(x)$}
        \label{fig:f}
    \end{subfigure}
    \hfill
    \begin{subfigure}[b]{0.3\textwidth}
        \centering
        \includegraphics[width=\textwidth]{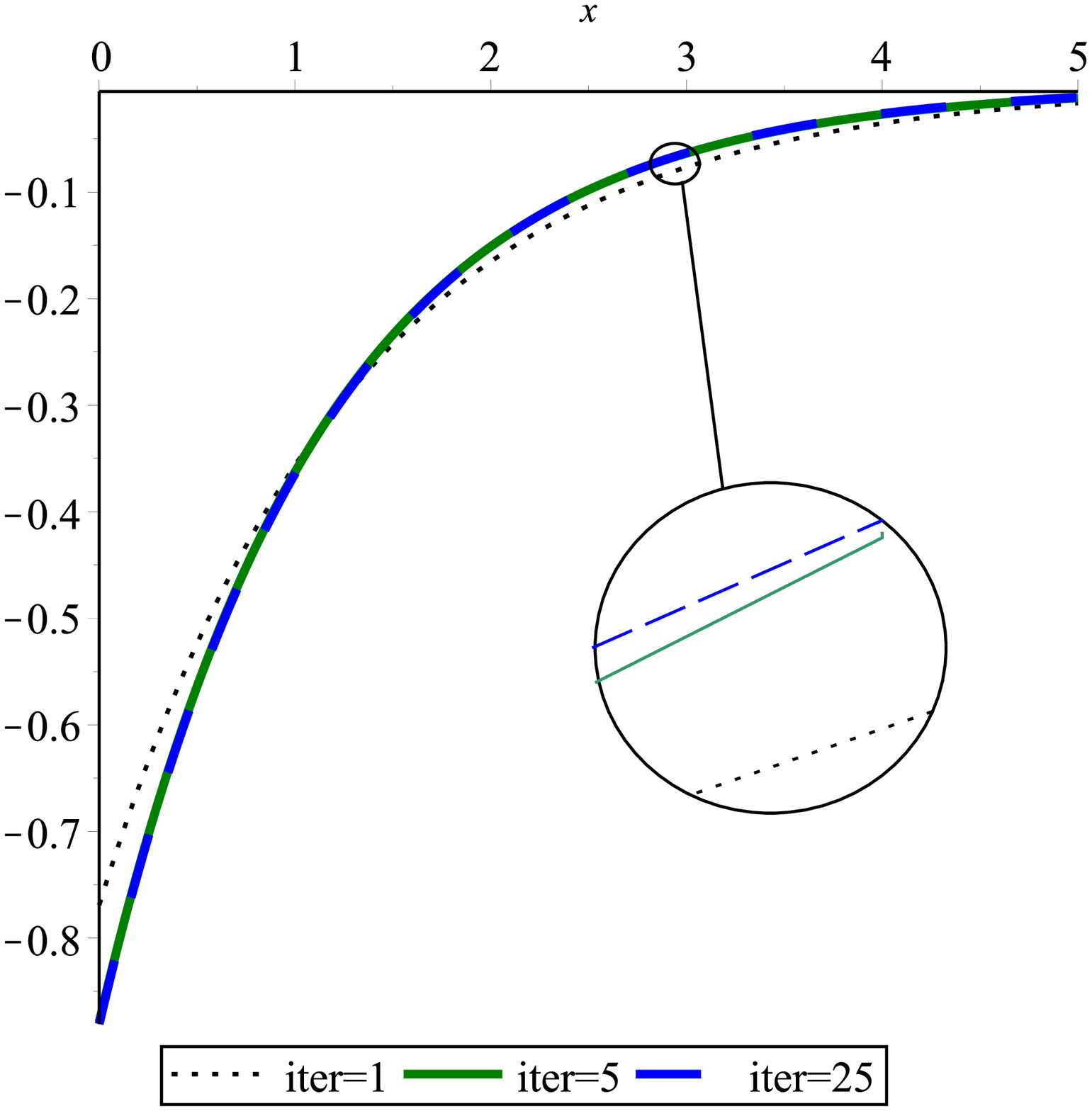}
        \caption{$f''(x)$}
        \label{fig:g}
    \end{subfigure}
     \hfill
    \begin{subfigure}[b]{0.4\textwidth}
        \centering
        \includegraphics[width=\textwidth]{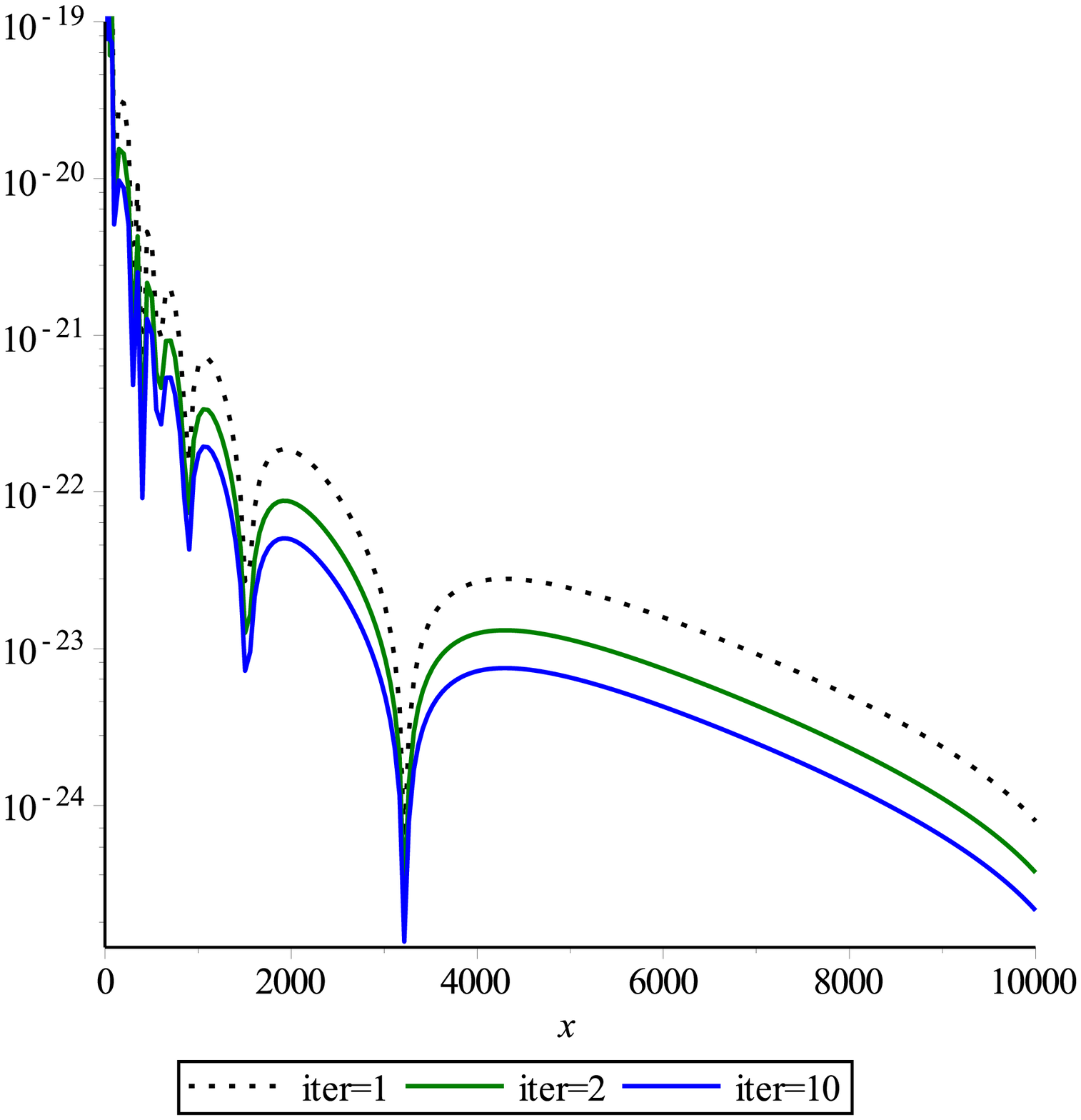}
        \caption{Logarithm of the absolute value of residual.}
        \label{fig:h}
    \end{subfigure}
    \caption{Diagram of $f(x)$, $f'(x)$, $f''(x)$ and $Log(|Res(x)|)$ for $N=50$ in different iterations. $\delta=0.1$, $\varepsilon=0.3$, $\alpha=1,\beta=1$, $L=15$. }
    \label{fig:three graphs}
\end{figure}

\begin{figure}[H]
    \centering
    \begin{subfigure}[b]{0.45\textwidth}
        \centering
        \includegraphics[width=\textwidth]{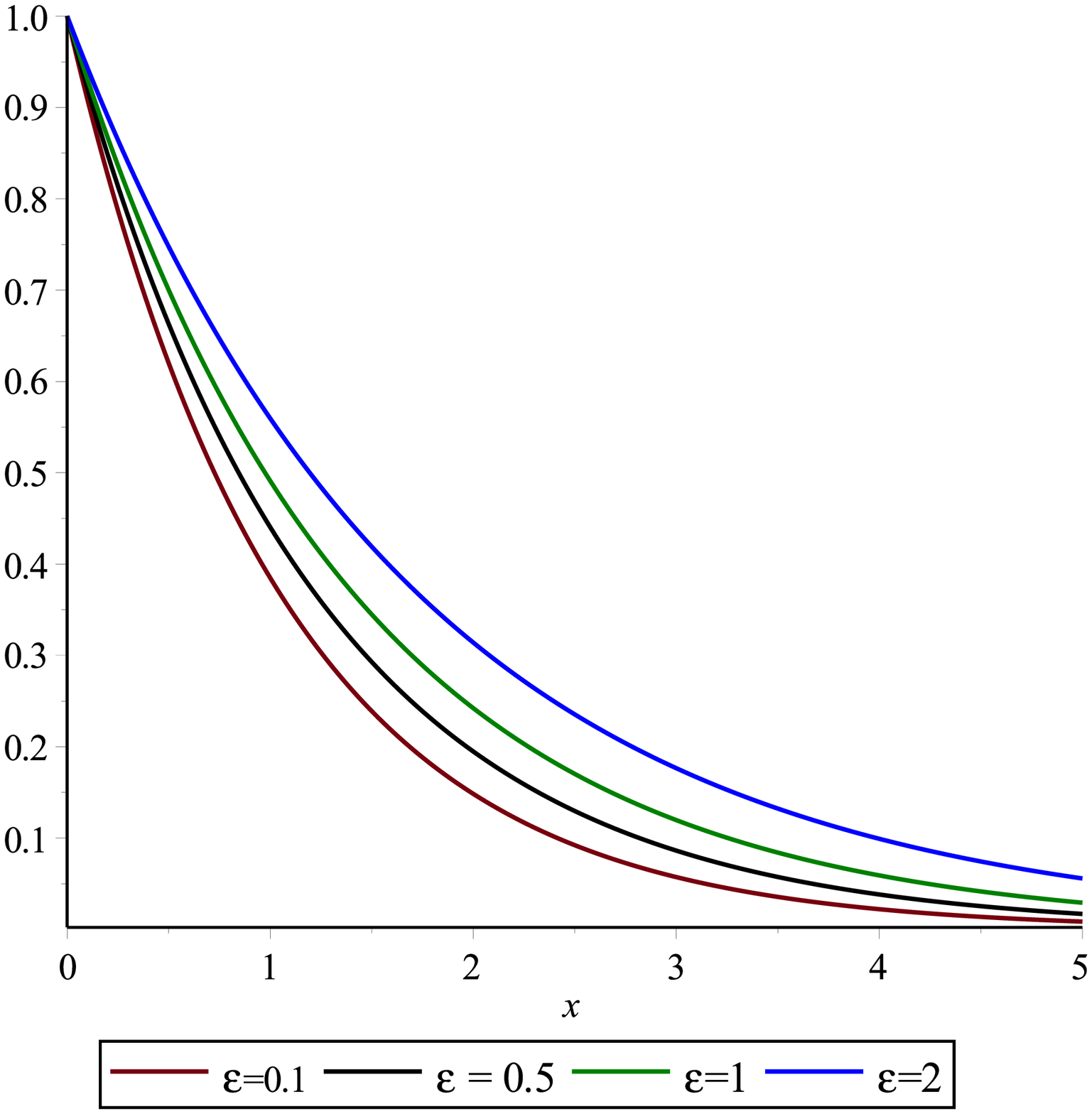}
        \caption{$f'(x)$}
        \label{fig:i}
    \end{subfigure}
    \hfill
    \begin{subfigure}[b]{0.45\textwidth}
        \centering
        \includegraphics[width=\textwidth]{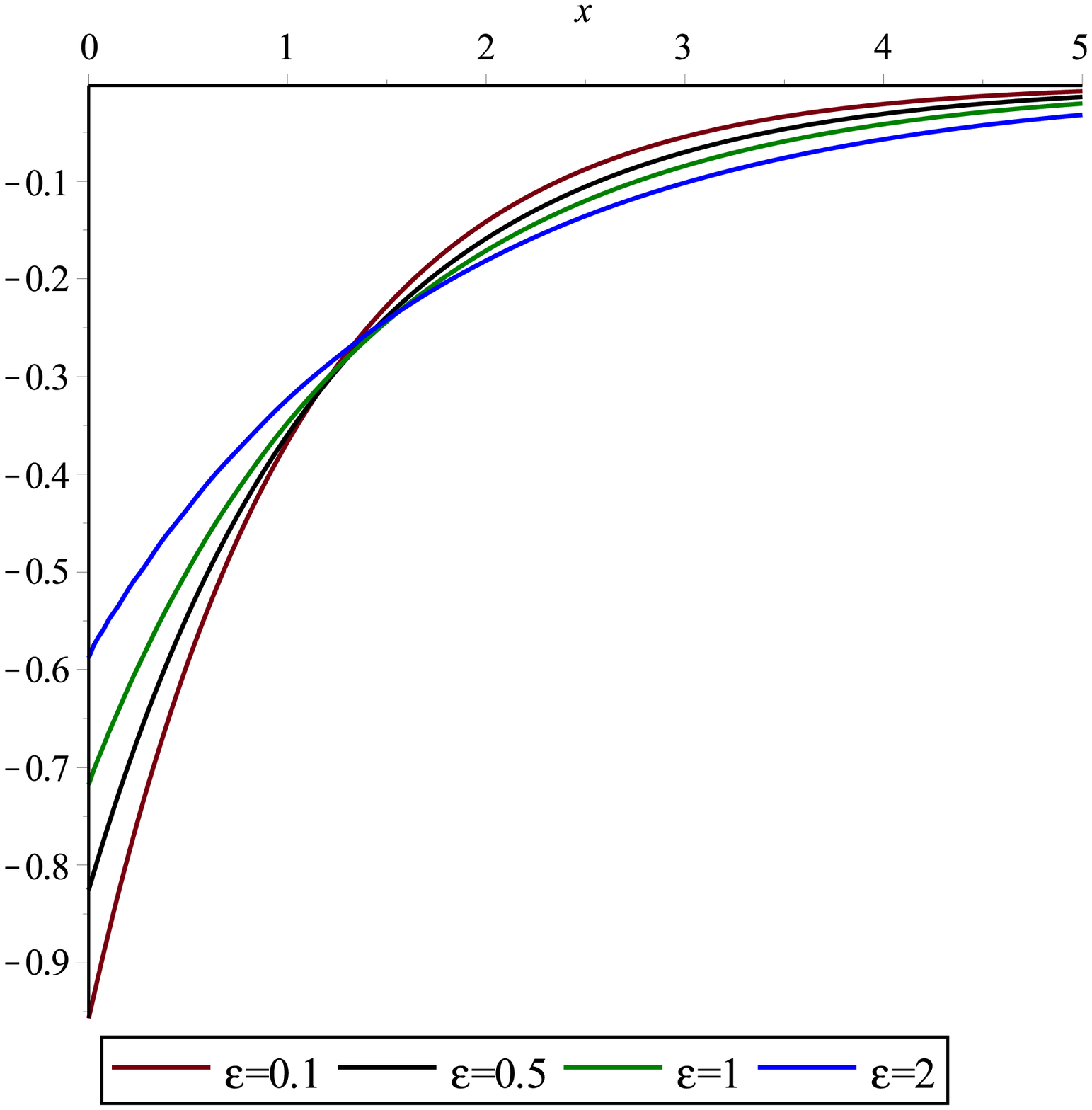}
        \caption{$f''(x)$}
        \label{fig:j}
    \end{subfigure}
    \caption{ Diagram of the effect of $\varepsilon$ on $f'(x)$ and $f''(x)$  in $15$th iteration. $\delta=0.2$, $N=50$, $\delta=0.2$ ,$\alpha=1,\beta=1$, $L=15$. 
}
    \label{fig:k}
\end{figure}

After solving the system of equations mentioned in Eq.(\ref{equation:eq6}) we listed the values of resulted solutions. In Table(\ref{tab:numericalresult1}) the result of $f(x)$,$f'(x)$ and $f''(x)$ for different $N$(Collocation points) are reported.  \\
As discussed earlier, the presented error is implemented without truncating domain. In order to show the accuracy of the presented mthod in Table (\ref{tab:compar}) we made a comparison  with other methods: IRBFs-QLM method \cite{lotfibaba}; Rational Chebyshev collocation method with truncation in domain \cite{parandmoayeri}; Runge-Kutta method; Collocation method with basis functions $\{1,x,x^2,...\}$ and truncation in domain \cite{lotfibaba}. 
From this table can draw a conclusion that since
the Runge-Kutta 4th order method is  exact at least up to 4 digits\cite{lotfibaba}, the presented method with no truncation is more accurate than Rational Chebychev collocation method endowed with truncation in domain. This is also   more accurate than CM (collocation method with  basis functions $\{1,x,x^2,...\}$).\\
In Fig.(\ref{fig:a},\ref{fig:b},\ref{fig:c} and \ref{fig:d}) and fig.(\ref{fig:e},\ref{fig:f},\ref{fig:g} and \ref{fig:h}) we showed that as the steps of iteration in QLM or Collocation points ($N$) increases the obtained result is better, the diagrams are smoother and the Logarithm of the absolute residual is closer to zero.\\ 
In Fig.(\ref{fig:i}) and  fig.(\ref{fig:j})  the effects of the non-Newtonian fluid parameters i.e $\delta$ and $\varepsilon$ on the non-dimensional velocity and stress field ($f'(x) $ and $f''(x)$) are presented. It is concluded that when $\varepsilon$ increases or $\delta$ decreases, $f'(x)$ increases. 
if $\delta$ increases or $\varepsilon$  decreases, as a result $f''(x)$ will increase.\\
If the obtained solution ($f_{i+1}(x)$) is put in the original equation (Eq. (\ref{equation:eq3})), the residual function will be defined as
\begin{equation}
res(x)=f_{i+1}(x)f_{i+1}''(x)+(1+\varepsilon)f_{i+1}'''(x)-\varepsilon\delta f_{i+1}''(x)^2f_{i+1}'''(x)-f_{i+1}'(x)^2,
\end{equation}   
The correctness and convergence of the method is illustrated by plotting residual function ($res(x)$) and the values of coefficients existed in  Eq.(\ref{equation:eq5}) are also denoted.
To calculate an approximation of the optimal value of parameter L we used the following fact:  \\
"The experimental trial-and-error method(Optimizing  Infinite  Interval Map Parameter)  \cite{boyd2001}: Plot  the  coefficients $a_j$ versus  degree  on  a  log-linear  plot.   If  the  graph abruptly flattens at some $N$, then this implies that $L$ is TOO SMALL for the  given $N$,  and  one  should  increase $L$ until  the  flattening  is  postponed to $j=N$."

\begin{center}
\begin{figure}[H]
\centering
 
 \includegraphics[scale=0.23]{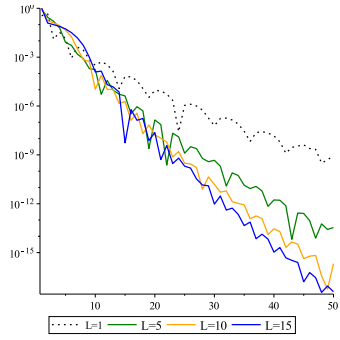}
  \caption{ Logarithm of the absolute value of coefficients  $\big(Log(|a_j^{i+1}|),j=0,...,N\big)$ that show  the optimal value of $L$ in the $25$th iteration. $N=50$, $\varepsilon=0.3$, $\delta=0.1$, $\alpha=1, \beta=1$. }
  \label{fig:optimaL}
\end{figure}

\end{center}

Based on this fact, the optimal parameter L is shown in Fig. \ref{fig:optimaL}. This optimal value is seemed to be $L=15$. Additionally, because the diagram of the coefficient $a_j, j=0,...,N$  decreases, it is obvious that the presented  method is  convergent.
\begin{table}[H]
\begin{center}\caption{ The numerical results in the 15th iteration for $f$,$f'$,$f''$ for different $N$. $\delta=0.1$, $\varepsilon=0.3$, $\alpha=1,\beta=1$, $L=15$.   }\label{tab:numericalresult1}
\scriptsize\begin{tabular}{ |c|c c c c c|} 
\hline
&x&N=10&N=15&N=25&N=50\\
\hline
\multirow{11}{*}{f}&0.0&0.0000000000&0.0000000000&0.0000000000&0.0000000000\\
                  &0.5&0.4036738856&0.4045266250&0.4045235026&0.4045235024\\
                  &1  &0.6658626593&0.6652631608&0.6652578446&0.6652578448\\
                  &1.5&0.8353708975&0.8334184925&0.8334154070&0.8334154068\\
                  &2.0&0.9443573988&0.9418982540&0.9418951684&0.9418951684\\
                  &2.5&1.0144374602&1.0118879495&1.0118841902&1.0118841903\\
                  &3.0&1.0595784826&1.0570456452&1.0570418968&1.0570418967\\
                  &3.5&1.0887030794&1.0861819735&1.0861787456&1.0861787455\\
                  &4.0&1.1075132470&1.1049814404&1.1049787070&1.1049787070\\
                  &4.5&1.1196651637&1.1171115798&1.1171090443&1.1171090444\\
                  &5.0&1.1275118485&1.1249385493&1.1249359382&1.1249359382\\
                  
\hline
\multirow{11}{*}{f'} &0.0 &1.0000000000&1.0000000000&1.0000000000&1.0000000000\\
                     &0.5 &0.6450803420&0.6442642937&0.6442495479&0.6442495461\\
                     &1.0 &0.4189562550&0.4154133505&0.4154176707&0.4154176711\\
                     &1.5 &0.2697139887&0.2679636516&0.2679659918&0.2679659914\\
                     &2.0 &0.1733336877&0.1728814969&0.1728800369&0.1728800375\\
                     &2.5 &0.1115538147&0.1115431922&0.1115424083&0.1115424080\\
                     &3.0 &0.0719247489&0.0719686778&0.0719693974&0.0719693972\\
                     &3.5 &0.0464362346&0.0464355030&0.0464366676&0.0464366678\\
                     &4.0 &0.0299991498&0.0299616664&0.0299623908&0.0299623909\\
                     &4.5 &0.0193776822&0.0193326322&0.0193327138&0.0193327137\\
                     &5.0 &0.0125058088&0.0124744388&0.0124741099&0.0124741098\\

\hline
\multirow{11}{*}{f''}  &0.0 &-0.9192031329&-0.8809437796&-0.8807848835&-0.8807849724\\
                       &0.5 &-0.5513725007&-0.5658276136&-0.5658037464&-0.5658037300\\
                       &1.0 &-0.3662977979&-0.3643882166&-0.3643723455&-0.3643723552\\
                       &1.5 &-0.2385329622&-0.2348981506&-0.2349101459&-0.2349101417\\
                       &2.5 &-0.0981133286&-0.0977533531&-0.0977500466&-0.0977500475\\
                       &3.0 &-0.0630210169&-0.0630697634&-0.0630676466&-0.0630676457\\
                       &3.5 &-0.0405944573&-0.0406920734&-0.0406922931&-0.0406922928\\
                       &4.0 &-0.0262103556&-0.0262544563&-0.0262557467&-0.0262557472\\
                       &4.5 &-0.0169503305&-0.0169398783&-0.0169410136&-0.0169410139\\
                       &5.0 &-0.0109702242&-0.0109304115&-0.0109308917&-0.0109308915\\

\hline

\end{tabular}
\end{center}
\end{table}

\begin{table}[H]
\begin{center}\caption{The value of $f'$ by presented method and other methods in the 15th iteration. $N=10$, $\delta=0.1$, $\varepsilon=0.3$, $\alpha=1,\beta=1$ and $L=15$.}\label{tab:compar}
\scriptsize\begin{tabular}{ c|c c c c c} 

x&	CM \cite{lotfibaba}&	Runge-katta \cite{lotfibaba} &	IRBFs-QLM \cite{lotfibaba}&	Rational Chebyshev \cite{parandmoayeri} &Presented method\\
\hline
0&	0.99999	&1	&1	&1&1\\
1&	0.41411	&0.41542	&0.41542&	0.4158220240&0.4189562550\\
2&	0.16994	&0.17288	&0.17288&	0.1737538392
&0.1733336877\\
3&	0.06748	&0.07197	&0.07197&	0.0732981252
&0.0719247489\\
4&	0.02426	&0.02996	&0.02996&	0.0316552721
&0.0299991498\\
5&	0.00593	&0.01247	&0.01247&0.0144254667&0.0125058088\\
\end{tabular}
\end{center}
\end{table}

\begin{figure}[H]
    \centering
    \begin{subfigure}[b]{0.45\textwidth}
        \centering
        \includegraphics[width=\textwidth]{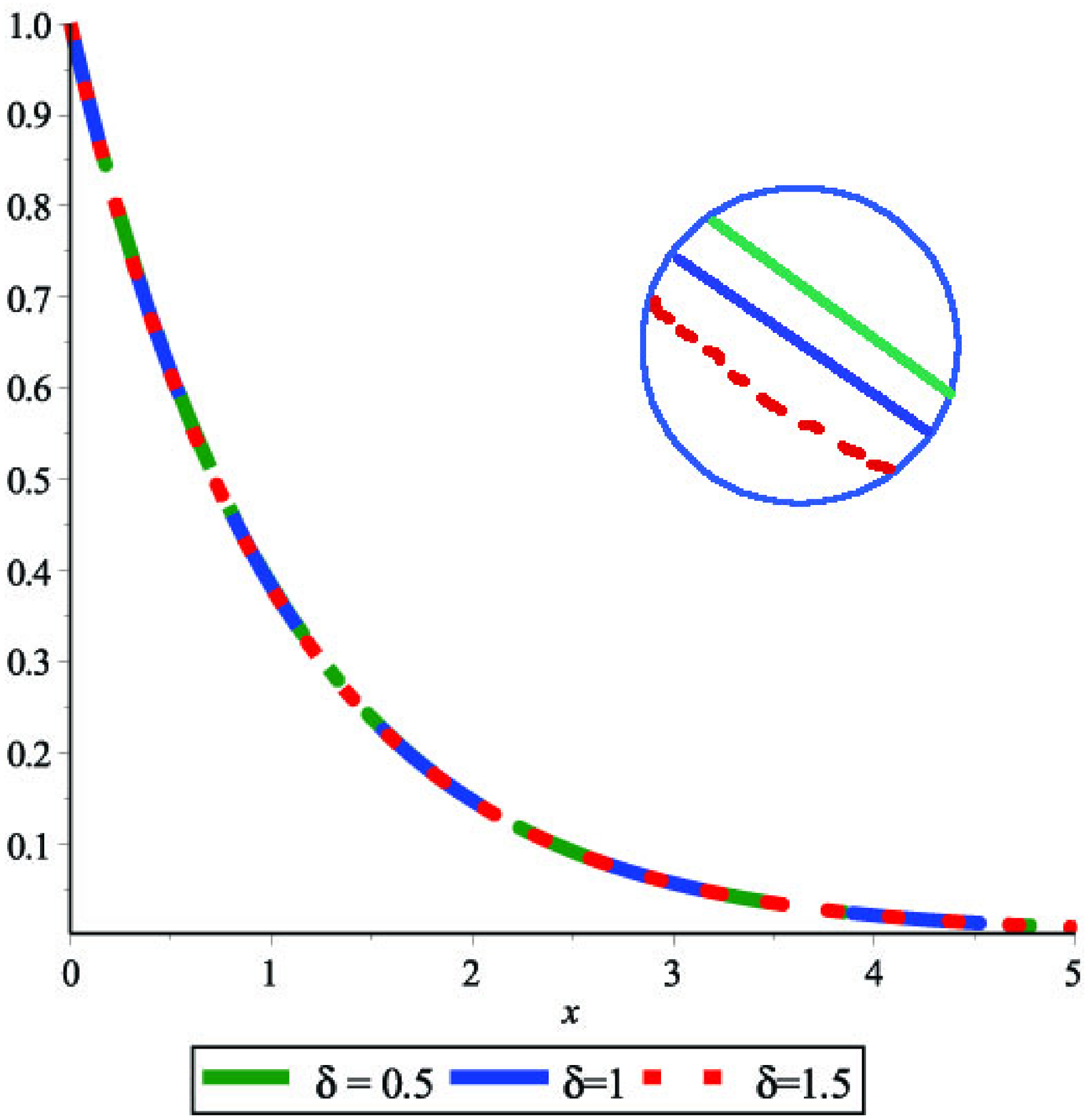}
        \caption{$f'(x)$}
        \label{fig:l}
    \end{subfigure}
    \hfill
    \begin{subfigure}[b]{0.45\textwidth}
        \centering
        \includegraphics[width=\textwidth]{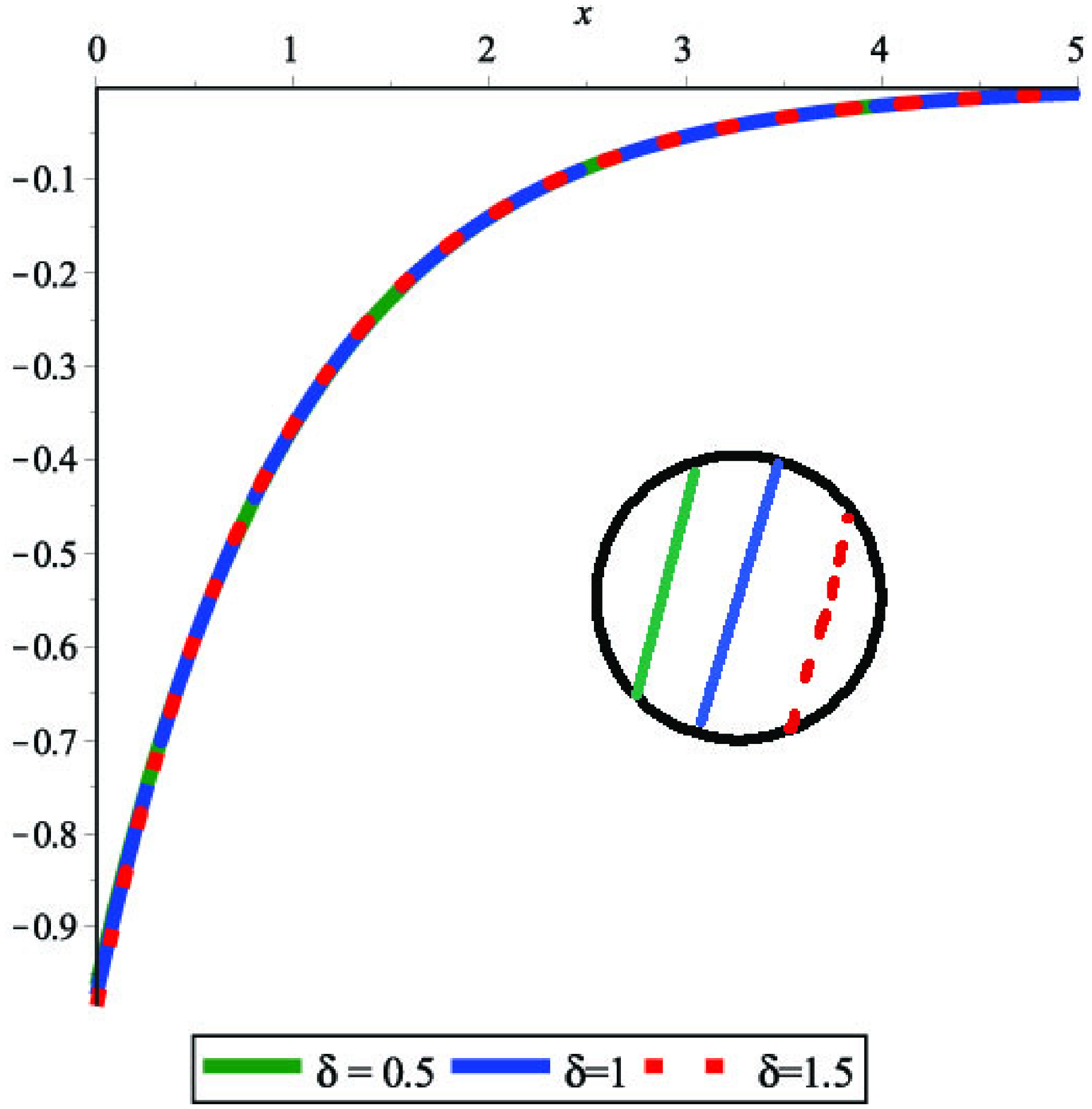}
        \caption{$f''(x)$}
        \label{fig:m}
    \end{subfigure}
    \caption{ Diagram of the effect of $\delta$ on $f'(x)$ and $f''(x)$  in the $15$th iteration. $N=50$, $\varepsilon=0.1$ ,$\alpha=1,\beta=1$, $L=15$. 
}
    \label{fig:n}
\end{figure}

\begin{center}
\begin{figure}[H]
\centering
 
 \includegraphics[scale=0.23]{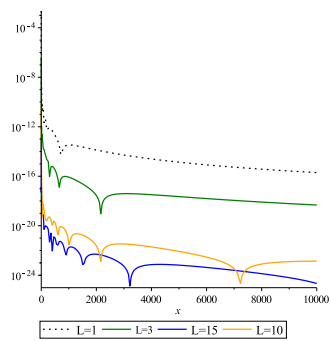}\label{fig:o}
  \caption{Logarithm of the absolute value of residual in the $15$th iteration by different $L$. $N=50$, $\varepsilon=0.3$, $\delta=0.1$, $\alpha=1,\beta=1$. }
\end{figure}
\end{center}

\section{conclusion}

In this paper, the RJ collocation method  applied for a nonlinear equation of momentum with an infinite condition  and the velocity and stress profiles of a non-Newtonian Eyring-Powell fluid over a linear stretching sheet is studied.
Due to encountering with an unbounded domain, we changed Jacobi polynomials and created RJ basis. Instead of truncation on domain that these days are used \cite{parandmoayeri}, we make the expansion  as if it satisfies all the conditions. By this expansion, the value of $f'(\infty)=0$ is reached implicitly. This presented approach has general applicability and can widely be used for solving the problems defined on an unbounded domain. We also used a method to linearize the nonlinear equation in question. As the measurement of $L$ is a crucial matter, its optimal value is also presented. To exhibit the convergence of the presented method, the effect and behavior of the coefficients are explained and this concluded that this method is convergent.   
Via graphical illustrations, The effect of different parameters on the
velocity profile is presented. It is shown that the velocity decreases by increasing the fluid material parameter ($\delta$) and increases by increasing the Eyring-Powell fluid material parameter ($\epsilon$).

\section*{Appendix}

One can read Eq. (\ref{equation:eq3}) as
\begin{equation*}
f'''(x)=\frac{f''(x)^2-f(x)f''(x)}{(1-\varepsilon )\delta\varepsilon f''(x)^2},
\end{equation*}
applying QLM technique, discussed in section \ref{section:QLM}, we obtain the approximation of Eq. (\ref{equation:eq3}), at step $i+1$, as like as
\begin{eqnarray}\label{equation:eq7}
&f'''_{i+1}(x)\approx\frac{f''_i(x)^2-f_i(x)f_i''(x)}{(1-\varepsilon )\delta\varepsilon f_i''(x)^2}+(f_{i+1}-f_i)\bigg(\frac{-f_i''(x)}{(1-\delta)-\delta\varepsilon f_i''(x)^2}\bigg)\nonumber\\
&+(f_{i+1}'(x)-f_{i}'(x))\bigg(\frac{2f_i'(x)}{(1-\delta)-\delta\varepsilon f_i''(x)^2}\bigg)\nonumber\\
&+(f_{i+1}''(x)-f_{i}''(x))\bigg(\frac{-f_i(x)-(f_i'(x)^2-f_if_i''(x))((1-\delta)-\delta\varepsilon f_i''(x)^2)}{((1-\delta)-\delta\varepsilon f_i''(x)^2)^2}\bigg),
\end{eqnarray}  
based on the philosophy of QLM technique, $f_i(x)$ is known, and  $f_{i+1}$ is being attempted to be solved. By factorization, the simpler form of the Eq. (\ref{equation:eq7}) can be shown as: 
\begin{equation*}
  g_1(.)+g_2(.)f_{i+1}(x)+g_3(.)f'_{i+1}(x)+g_4(.)f''_{i+1}(x)-f'''_{i+1}(x)\approx 0,
  \end{equation*}
where 
\begin{eqnarray*}
&g_1(.)=\frac{1}{-\varepsilon\delta f_i''(x)^2-\varepsilon+1}\bigg(-f_i(x)f_i''(x)-f_i'(x)^2-\\
&\frac{x}{x^2+1}f_i(x)''+2\frac{-x^2+1}{(x^2+1)^2}-\frac{2x^3-6x}{(x^2+1)^3}f_i(x)\bigg),
\end{eqnarray*}

\begin{equation*}
g_2(.)=\frac{1}{-\varepsilon\delta f_i''(x)^2-\varepsilon+1}\bigg((\frac{-x^2}{x^2+1})f_i''(x)+\frac{4x}{(x^2+1)^2}f_i'(x)-(\frac{-6x^2+2}{(x^2+1)3})f_i(x)''\bigg),
\end{equation*}

\begin{equation*}
g_3(.)=\frac{1}{-\varepsilon\delta f_i''(x)^2-\varepsilon+1}\bigg((\frac{-2x^4-2x^2}{(x^2+1)^2})f_i'(x)+(\frac{-4x^3-4x}{(x^2+1)^3})f_i(x)\bigg),
\end{equation*}

\begin{equation*}
g_4(.)=\frac{1}{-\varepsilon\delta f_i''(x)^2-\varepsilon+1}\bigg((\frac{-x^6-2x^4-x^2}{(x^2+1)^3})f_i(x)\bigg).
\end{equation*}

for simplicity [(.) $\equiv$ $\big(x,f_i(x),f'_i(x),f''_i(x)\big)$].


\bibliographystyle{spmpsci}      


\end{document}